\documentclass[11pt]{amsart}

\usepackage{amsmath,amssymb,amsthm}
\usepackage[margin=1.1in]{geometry}
\usepackage{hyperref}

\theoremstyle{plain}
\newtheorem{theorem}{Theorem}[section]
\newtheorem{lemma}[theorem]{Lemma}
\newtheorem{proposition}[theorem]{Proposition}
\newtheorem{corollary}[theorem]{Corollary}
\theoremstyle{remark}
\newtheorem*{remark}{Remark}

\newcommand{\OP}{\mathit{OP}}
\newcommand{\OR}{\mathit{OR}}
\newcommand{\GR}{\mathcal{R}}
\newcommand{\GL}{\mathcal{L}}
\newcommand{\GH}{\mathcal{H}}
\newcommand{\GD}{\mathcal{D}}
\newcommand{\GJ}{\mathcal{J}}
\newcommand{\GG}{\mathcal{G}}
\newcommand{\F}{\mathrm{F}}
\newcommand{\gen}[1]{\langle #1\rangle}
\newcommand{\tmap}[6]{\begin{pmatrix} #1 & #2 & #3\\ #4 & #5 & #6\end{pmatrix}}

\title{The Monoid of Orientation-preserving Mappings on a Chain}

\author{P.~M.~Catarino}
\address{Departamento de Matem\'atica, Universidade de Tr\'as-os-Montes e Alto Douro,
Quinta de Prados, 5000 Vila Real, Portugal}
\email{pcatarin@utad.pt}

\author{P.~M.~Higgins}
\address{Department of Mathematics, University of Essex, Wivenhoe Park,
Colchester CO4 3SQ, U.K.}
\email{peteh@essex.ac.uk}

\thanks{This is a transcription of the article published in
\emph{Semigroup Forum} \textbf{58} (1999), 190--206
(received May 13, 1996; communicated by John M.~Howie).
Minor typographical corrections have been made in this version.}


\begin{document}

\begin{abstract}
We introduce new kinds of regular monoids of transformations consisting of
mappings which preserve or reverse the orientation of a finite cycle.
\end{abstract}

\maketitle

% ===============================================================
\section{Introduction and preliminary concepts}
% ===============================================================

For standard terms in semigroup theory see, for example, \cite{CP},
\cite{H92}, \cite{Howie95}. Let $X$ denote the set $\{1,2,\dots,n\}$
considered with the standard ordering. Let $T_n$ denote the full
transformation semigroup on $X$, that is the semigroup of all mappings
$\alpha\colon X\longrightarrow X$ under composition. We shall denote
composition simply by juxtaposition and shall adopt a right mapping
convention; thus if $\alpha,\beta\in T_n$ the symbol $\alpha\beta$ denotes
the mapping in $T_n$ obtained by performing first $\alpha$ and then
$\beta$. The range and the set of fixed points of $\alpha\in T_n$ will be
denoted by $X\alpha$ and $\F(\alpha)$, respectively. The order of
$X\alpha$ will be called the \emph{rank} of $\alpha$. We shall denote by
$S_n$ the symmetric group on $X$. In this paper we frequently use the
cyclic group $Z_n$ with generator $a=(12\cdots n)$. We shall denote by
$C(n,q)$ the number of $q$-combinations of an $n$-set. We shall denote by
$Y\cdot Z$ the set product $\{yz : y\in Y, z\in Z\}$. We shall call
$\alpha$ an \emph{order-preserving} mapping if $i\le j$ implies
$i\alpha\le j\alpha$ $(i,j\in X)$ and let $O_n$ denote the submonoid of
$T_n$ consisting of all order-preserving mappings on $X$. The subsemigroup
generated by a non-empty subset $T$ of a semigroup $S$ is denoted by
$\gen{T}$. We will denote by $E=E(S)$ the set of idempotents of $S$. The
set of inverses of $y\in S$ is denoted by $V(y)$.

Various special subsemigroups of $T_n$ have been studied by many authors;
see for example, \cite{H93}, \cite{H94}, \cite{H95}, \cite{Howie66},
\cite{Howie71}, \cite{HS}, \cite{Umar}. The purpose of this paper is to
introduce a new subsemigroup of $T_n$ containing $O_n$ which will be
denoted by $\OP_n$ and its elements will be called
\emph{orientation-preserving} mappings. Also, we introduce a subsemigroup
$P_n=\OP_n\cup\OR_n$ of $T_n$, where $\OR_n$ denotes the collection of all
\emph{orientation-reversing} mappings.

In Section~2 we introduce the monoid $\OP_n$ and show that any element of
$\OP_n$ can be written as a product of a power of $a$ followed by an
order-preserving mapping. We also show that every kernel class of
$\alpha\in\OP_n$ is convex. In Section~3 we show that $\OP_n$ is a regular
semigroup, and we also characterise Green's relations on $\OP_n$. In
Section~4 we consider the subsemigroups of $\OP_n$ conjugate with respect
to $Z_n$. From this we obtain a description of the core of $\OP_n$. In
Section~5 we show that $P_n$ is a regular semigroup and also give some its
properties derived from the results about $\OP_n$.

First we need some definitions and results concerning sequences of
elements from a chain. In this paper it will always be clear from context
when additions are taken modulo $n$ (or modulo $t$ where $t$ is the number
of elements of any sequence). A (finite) sequence from a chain
$A=(a_1,a_2,\dots,a_t)$ is \emph{cyclic} or has \emph{clockwise
orientation} if there exists no more than one subscript $i$ such that
$a_i>a_{i+1}$. We say that $A=(a_1,a_2,\dots,a_t)$ is anti-cyclic or has
anti-clockwise orientation if there exists no more than one subscript $i$
such that $a_i<a_{i+1}$. For example $(1,1,2,5,1,1,1)$ is cyclic,
$(1,1,1,5,2,1,1)$ is anti-cyclic whilst $(1,1,1,5,1,1,1)$ is both cyclic
and anti-cyclic.

For any result concerning cyclic sequences there is a dual result for
anti-cyclic sequences obtained by reversing the direction of all
inequalities. We shall not prove these results in the sequel but consider
them to be implicitly proved once the results on cyclic sequences are
established.

\begin{remark}
Note that a sequence $A=(a_1,a_2,\dots,a_t)$ is cyclic if and only if
there exists $i\in\{0,\dots,t-1\}$ such that
\[
a_{i+1}\le\cdots\le a_t\le a_1\le\cdots\le a_i
\]
and $i$ is unique unless the sequence is constant (where we adopt the
convention that $\le a_1\le\cdots\le a_i$ is vacuous if $i=0$).
\end{remark}

The following lemma is immediate by inspection using the remark above:

\begin{lemma}
Let $A$ be any cyclic \textup{[}anti-cyclic\textup{]} sequence. Then $A$
is anti-cyclic \textup{[}cyclic\textup{]} if and only if $A$ has no more
than two distinct values.
\end{lemma}

If $A=(a_1,a_2,\dots,a_t)$ is any sequence then we denote by $A^r$ the
sequence $(a_t,a_{t-1},\dots,a_1)$, called the \emph{reversed sequence} of
$A$.

We have three more results about cyclic [anti-cyclic] sequences, the
first two of which are obvious.

\begin{lemma}
Let $A=(a_1,a_2,\dots,a_t)$ be any sequence from $X$. Then $A$ is cyclic
\textup{[}anti-cyclic\textup{]} if and only if $A^r$ is anti-cyclic
\textup{[}cyclic\textup{]}.
\end{lemma}

\begin{lemma}
If $(a_1,a_2,\dots,a_t)$ is cyclic \textup{[}anti-cyclic\textup{]} then
so is
\begin{enumerate}
\item[(a)] the subsequence $(a_{i_1},a_{i_2},\dots,a_{i_r})$
$(i_1<i_2<\cdots<i_r)$;
\item[(b)] $(a_j,a_{j+1},\dots,a_t,a_1,\dots,a_{j-1})$ for all
$1\le j\le t$.
\end{enumerate}
\end{lemma}

Let $k\in\{0,1,\dots,n-1\}$. Define a total order $\le_k$ on $X$ by
\[
k+1\le_k k+2\le_k\cdots\le_k n\le_k 1\le_k\cdots\le_k k
\]
(where it is understood that $\le_k 1\le_k\cdots\le_k k$ is vacuous if
$k=0$). We write $i<_k j$ if $i\le_k j$ and $i\ne j$. Note that
$i\le_k j$ if and only if $i-k\le j-k$.

We say that $A=(a_1,a_2,\dots,a_t)$ is \emph{cyclic}
[\emph{anti-cyclic}] \emph{with respect to} $\le_k$ if there exists no
more than one subscript $i$ such that $a_{i+1}<_k a_i$
[$a_i<_k a_{i+1}$]. Clearly $A$ is cyclic with respect to $\le_k$ if and
only if there exists $j\in\{0,\dots,t-1\}$ such that
$a_{j+1}\le_k\cdots\le_k a_t\le_k a_1\le_k\cdots\le_k a_j$.

\begin{lemma}\label{lem:1.4}
Let $A=(a_1,a_2,\dots,a_t)$ be any sequence of elements from $X$. Then the
following are equivalent:
\begin{enumerate}
\item[(a)] $A$ is cyclic with respect to $\le_0\;=\;\le$.
\item[(b)] $A$ is cyclic with respect to $\le_k$ for some $k$.
\item[(c)] $A$ is cyclic with respect to $\le_k$ for all $k$.
\end{enumerate}
\end{lemma}

\begin{proof}
Clearly we only need to prove that (b) implies (c). Suppose that (b) holds
so that $A$ is cyclic with respect to $\le_k$ for some $k$. Then, for some
$j$,
\[
k+1\le_k a_{j+1}\le_k\cdots\le_k a_t\le_k a_1\le_k\cdots\le_k a_j\le_k k.
\]
Either $k+2\le_k a_{j+1}$, in which case
\[
a_{j+1}\le_{k+1}\cdots\le_{k+1}a_t\le_{k+1}a_1\le_{k+1}\cdots\le_{k+1}a_j,
\]
or $k+1=a_{j+1}=\cdots=a_{j+s}<_k k+2\le_k a_{j+s+1}\le_k\cdots\le_k a_j$
for some $s\in\{1,\dots,t\}$, in which case
\[
a_{j+s+1}\le_{k+1}\cdots\le_{k+1}a_j\le_{k+1}\cdots\le_{k+1}a_{j+s}.
\]
In either case $A$ is cyclic with respect to $\le_{k+1}$, and then (c)
follows by induction.
\end{proof}

% ===============================================================
\section{The monoid \texorpdfstring{$\OP_n$}{OPn}}
% ===============================================================

Let $\alpha\in T_n$. We say that $\alpha$ is an
\emph{orientation-preserving} mapping on $X$ if the sequence
$(1\alpha,2\alpha,\dots,n\alpha)$ is cyclic. We note, by
Lemma~\ref{lem:1.4}, that this sequence is then cyclic with respect to
$\le_k$ for all $0\le k\le n-1$. The collection of all
orientation-preserving mappings on $X$ will be denoted by $\OP_n$.

\begin{lemma}\label{lem:2.1}
Let $\alpha\in\OP_n$ and let $(b_1,\dots,b_t)$ be a cyclic sequence of
members of $X$. Then the sequence $(b_1\alpha,\dots,b_t\alpha)$ is also
cyclic.
\end{lemma}

\begin{proof}
We have, for some $i\in\{0,\dots,t-1\}$
\[
b_{i+1}\le b_{i+2}\le\cdots\le b_t\le b_1\le\cdots\le b_i,
\]
and, for some $j\in\{0,\dots,n-1\}$,
\[
(j+1)\alpha\le(j+2)\alpha\le\cdots\le n\alpha\le 1\alpha\le\cdots\le j\alpha.
\]
In particular $\alpha$ preserves $\le$ on $\{j+1,\dots,n\}$ and on
$\{1,\dots,j\}$. Thus if $b_{i+1}\in\{j+1,\dots,n\}$ or
$b_i\in\{1,\dots,j\}$ then certainly
\[
b_{i+1}\alpha\le b_{i+2}\alpha\le\cdots\le b_t\alpha\le b_1\alpha\le\cdots\le b_i\alpha
\]
and we have finished. Suppose then that $b_{i+1}<j+1$ and $b_i>j$. Then,
for some $r$, either
\[
i+1<r\le t,\quad b_{i+1},\dots,b_{r-1}\in\{1,\dots,j\},
\]
or
\[
1\le r\le i,\quad b_{i+1},\dots,b_t,b_1,\dots,b_{r-1}\in\{1,\dots,j\},\quad
b_r,\dots,b_i\in\{j+1,\dots,n\}.
\]
In either case
\[
b_r\alpha\le\cdots\le b_t\alpha\le b_1\alpha\le\cdots\le b_{r-1}\alpha,
\]
as required.
\end{proof}

As a corollary we immediately obtain the following result:

\begin{theorem}
$\OP_n$ is a submonoid of $T_n$ containing $O_n$.
\end{theorem}

The following is immediate from the definitions of order-preserving
mapping on $X$ and cyclic sequence:

\begin{proposition}\label{prop:2.3}
For non-constant $\alpha\in\OP_n$, $\alpha$ is an order-preserving mapping
if and only if $1\alpha<n\alpha$.
\end{proposition}

We may regard the members of $X$ as being placed clockwise around the
circumference of a circle so that the integer $i$ lies between $i-1$ and
$i+1$ (reduced modulo $n$). Any sequence of 3 distinct members, $i,j,k$,
of $X$ then defines either a clockwise orientation or an anti-clockwise
orientation of the circle according as $(i,j,k)$ is cyclic or anti-cyclic.
The following lemma give a relation between an orientation-preserving
mapping and triples of 3 distinct members of $X$, and follows immediately
from Lemmas~2.1 and~1.2,

\begin{proposition}
Let $\alpha\in\OP_n$ and $(i,j,k)$ be any triple of $3$ distinct members of
$X$. If the entries are distinct, then the triple
$(i\alpha,j\alpha,k\alpha)$ defines the same orientation as $(i,j,k)$.
\end{proposition}

An interval $[i,k]$ for $i,k\in X$ may be expressed \emph{clockwise} by
$[i,k]=\{i,i+1,\dots,k-1,k\}$ or \emph{anti-clockwise} by
$[i,k]=\{i,i-1,\dots,k+1,k\}$. Henceforth, all intervals are expressed
clockwise. A set $K\subseteq X$ is \emph{convex} if $K$ has the form
$[i,i+t]$, for some $i\in X$ and $0\le t\le n-1$. If $\alpha\in\OP_n$, let
us call the subsets $x\alpha^{-1}$ $(x\in X\alpha)$ of $X$ the kernel
classes of $\alpha$; they are the classes of equivalence relation
\[
\Pi_\alpha=\{(x,y)\in X\times X : x\alpha=y\alpha\}.
\]
Following Howie \cite{Howie71}, we shall say that a kernel class is
\emph{trivial} if it consists of a single element, and we use the notation
$\ker\alpha$ for the set of kernel classes of $\alpha$. It is easy to show
that every kernel class of $\alpha\in O_n$ is convex; indeed they are
intervals of $X$. Later, we prove more generally that if $\alpha\in\OP_n$,
then every kernel class of $\alpha$ is convex also.

The $n$-cycle $a=(12\cdots n)$ lies in $\OP_n$ as the sequence of images
of the ordered set $X$ under $a$, $(2,\dots,n,1)$, is cyclic.

We consider $S=\gen{O_n\cup\{a\}}$, and let $g\in O_n$. Clearly, if $g$
is constant, then $a^rg=g\in O_n$, for all $0\le r\le n-1$. We claim that,
if $g$ is not a constant mapping, then $a^rg\notin O_n$, for all
$1\le r\le n-1$. Suppose that $g$ is not constant. Since $g$ is an
order-preserving mapping it follows that $i=1g<ng=j$. Now $n-r<n-r+1$ but
$(n-r)a^rg=ng=j>i=1g=(n-r+1)a^rg$, showing that $a^rg\notin O_n$. Hence the
next result follows:

\begin{lemma}\label{lem:2.5}
Let $g\in O_n$. Then $a^rg\in O_n$, $(1\le r\le n-1)$, if and only if $g$
is constant. \qed
\end{lemma}

Next we give an interesting factorization of $\OP_n$ and we show that any
orientation-preserving mapping is a product of a power of $a$ followed by
an order-preserving mapping.

\begin{theorem}\label{thm:2.6}
Each $\alpha\in\OP_n$ has a factorization $\alpha=a^rf$
$(0\le r\le n-1,\ f\in O_n)$ which is unique unless $\alpha$ is constant.
In particular $\OP_n=S$.
\end{theorem}

\begin{proof}
Let $\alpha\in\OP_n$. We claim that $\alpha=a^rf$ for some $f\in O_n$ and
$0\le r\le n-1$. If $\alpha$ is a constant mapping then we take
$f=\alpha$ and $r$ arbitrary. Suppose that $\alpha$ is not a constant
mapping. Then there exists $i\in X$ such that $i\alpha>(i+1)\alpha$ and we
have that
\[
\min_{1\le j\le n}j\alpha=(i+1)\alpha\le\cdots\le n\alpha\le 1\alpha\le\cdots\le i\alpha=\max_{1\le j\le n}j\alpha.
\]
We define $f$ by $xf=(x+i)\alpha$ for $x\in X$. Then $f$ lies in $O_n$ as
$1f=(1+i)\alpha\le 2f=(2+i)\alpha\le\cdots\le nf=(n+i)\alpha=i\alpha$. Put
$r=n-i$. Then $j(a^rf)=(j+r)f=(j+r+i)\alpha=(j+n)\alpha=j\alpha$,
$1\le j\le n$. Next suppose $\alpha$ is not constant and
$\alpha=a^kf=a^lg$, $(0\le k,l\le n-1)$, $f,g\in O_n$. Then
$f=a^{l-k}g$, so that, by Lemma~\ref{lem:2.5}, $l-k=0$, since $g$ is not
constant, whence $l=k$ and $f=g$.
\end{proof}

Now, using the known cardinalities of $O_n$ (see for example
\cite{Howie71}) and $Z_n$, we obtain the following result:

\begin{corollary}\label{cor:2.7}
$|S|=|\OP_n|=nC(2n-1,n-1)-n(n-1)$.
\end{corollary}

\begin{proof}
By above,
\[
|S|=|Z_n|(|O_n|-n)+n=n|O_n|-n^2+n=nC(2n-1,n-1)-n(n-1),
\]
as required.
\end{proof}

Next we show that the kernel classes of orientation-preserving mappings
are necessarily convex.

\begin{lemma}[Convexity Lemma]\label{lem:2.8}
Let $\alpha\in\OP_n$ and $K$ be any kernel class of $\alpha$. Then $K$ is
convex.
\end{lemma}

\begin{proof}
Let $\alpha\in\OP_n$. Then, for some $i$,
\[
(i+1)\alpha\le\cdots\le n\alpha\le 1\alpha\le\cdots\le i\alpha.
\]
Thus if $j\alpha=k\alpha$ for some $j\le_i k$ then
$j\alpha=(j+1)\alpha=\cdots=(k-1)\alpha=k\alpha$, which is sufficient to
prove that the kernel classes are convex.
\end{proof}

In \cite{Howie71} Howie showed that $|E(O_n)|=F_{2n}$, the $2n$th
Fibonacci number, through solution of a certain difference equation. The
next lemma gives another proof of this.

\begin{lemma}\label{lem:2.9}
Let $F_n$ be the $n$th Fibonacci number. Then, for all $n\ge1$,
\[
\text{(i)}\ |E(O_n)|=F_{2n};\qquad
\text{(ii)}\ |E(\OP_n)|=n+\sum_{k=1}^{n-1}k^2F_{2(n-k)}.
\]
\end{lemma}

\begin{proof}
(i) Let $e_n=|E(O_n)|$. Let $e\in E(O_{n+1})$. We consider two cases:

$Xe\subseteq X$: In this case $e|_X\in E(O_n)$. Now $ne\le(n+1)e$, as
$e$ is order-preserving, and $(n+1)e\le n$ as $Xe\subseteq X$. Applying
the order-preserving mapping $e$ to the latter inequality gives
$(n+1)e^2\le ne$, and since $e^2=e$ this says $(n+1)e\le ne$.
Therefore $(n+1)e=ne$ and so
$|\{e\in E(O_{n+1}) : Xe\subseteq X\}|=e_n$.

$Xe\not\subseteq X$: Let $i$ be least such $ie=(n+1)$. Then $je=n+1$, for
all $j\ge i$. In general $je\le i-1$, for all $j\le i-1$ and so
$(i-1)e\le i-1$. Hence,
$|\{e\in E(O_{n+1}) : i \text{ least such that } ie=(n+1)\}|=e_{i-1}$.
Summing over all values of $i$ and adding the previous case yields the
recurrence
\[
e_{n+1}=e_n+(e_n+e_{n-1}+\cdots+e_1+1),
\]
from which the result follows by induction.

(ii) Now let $e_n=|E(\OP_n)|$. We prove that the numbers of idempotents
$e$ in $\OP_n$ with $|1e^{-1}|=k$, $(1\le k\le n-1)$ is $k^2F_{2(n-k)}$,
from which the result will follow. We build such an idempotent
$e\in E(\OP_n)$ as follows:

There are $k$ choices for the convex set $1e^{-1}$:
$\{n-k+2,n-k+3,\dots,1\}$, $\{n-k+3,n-k+4,\dots,1,2\}$, $\dots$,
$\{1,2,\dots,k\}$. Once $1e^{-1}$ is chosen, there are $k$ choices for
$1e$. Suppose that $1e^{-1}=\{t,t+1,\dots,n,1,2,\dots,s\}$ say,
$(s=k+t-n-1)$. Put $1e^{-1}=K$. Then $e|_{X\setminus K}$ is an idempotent
on the chain $C=\{s+1,s+2,\dots,t-1\}$, $|C|=n-k$. By (i) it follows that
$e$ can be extended from $e|_K$ to an idempotent in $\OP_n$, in
$F_{2(n-k)}$ ways. Multiplying the number of choices in each of the 3
stages of the construction of $e$ gives the number of possibilities as
$kkF_{2(n-k)}=k^2F_{2(n-k)}$, as required.
\end{proof}

\begin{theorem}
For all $n\ge1$, $|E(\OP_n)|=F_{2n-1}+F_{2n+1}-(n^2-n+2)$.
\end{theorem}

\begin{proof}
We use the fact that for all $n\ge0$,
$F_0+F_2+\cdots+F_{2n}=F_{2n+1}-1$ and
$F_1+F_3+\cdots+F_{2n+1}=F_{2n+2}$, each of which can be proved readily by
induction. Let $e_n=|E(\OP_n)|$. We seek to simplify the formula given by
the previous lemma. First note that
$\sum_{k=0}^{n-1}F_{2(n-k)}=\sum_{k=1}^{n}F_{2k}=F_{2n+1}-1$ and
\begin{align*}
\sum_{k=1}^{n-1}kF_{2(n-k)}
&=\sum_{j=1}^{n-1}\sum_{k=j}^{n-1}F_{2(n-k)}
 =\sum_{j=1}^{n-1}(F_2+F_4+\cdots+F_{2(n-j)})\\
&=\sum_{j=1}^{n-1}(F_{2(n-j)+1}-1)
 =F_3+F_5+\cdots+F_{2n-1}-(n-1)\\
&=F_{2n}-n.
\end{align*}
Also we have
\begin{align*}
e_{n+1}-e_n
&=(n+1)-n+\sum_{k=1}^{n}k^2F_{2(n+1-k)}-\sum_{k=1}^{n-1}k^2F_{2(n-k)}\\
&=1+\sum_{k=1}^{n}k^2F_{2(n-(k-1))}-\sum_{k=1}^{n-1}k^2F_{2(n-k)}\\
&=1+\sum_{k=0}^{n-1}(k+1)^2F_{2(n-k)}-\sum_{k=0}^{n-1}k^2F_{2(n-k)}\\
&=1+\sum_{k=0}^{n-1}F_{2(n-k)}+2\sum_{k=1}^{n-1}kF_{2(n-k)}\\
&=1+(F_{2n+1}-1)+2(F_{2n}-n)\\
&=F_{2n}+F_{2n+2}-2n
\end{align*}
and so
\begin{align*}
e_{n+1}&=e_n+(F_{2n}+F_{2n+2})-2n=\cdots=\\
&=e_1+(F_2+F_4)+(F_4+F_6)+\cdots+(F_{2n-2}+F_{2n})+(F_{2n}+F_{2n+2})\\
&\qquad-2(n+(n-1)+(n-2)+\cdots+1)\\
&=2(F_2+F_4+\cdots+F_{2n})+F_{2n+2}-2\bigl(\tfrac12 n(n+1)\bigr)\\
&=2(F_{2n+1}-1)+F_{2n+2}-n(n+1)\\
&=F_{2n+1}+F_{2n+3}-(n^2+n+2),
\end{align*}
or $e_n=F_{2n-1}+F_{2n+1}-(n^2-n+2)$, as required.
\end{proof}

Next we describe $O_n\cdot Z_n$ and show that $O_n\cdot Z_n\ne\OP_n$.

\begin{proposition}
$O_n\cdot Z_n=\{\alpha\in\OP_n : 1\alpha\ne n\alpha\}\cup\{\alpha\in\OP_n : |X\alpha|=1\}$.
\end{proposition}

\begin{proof}
Let $\alpha\in O_n\cdot Z_n$ be such that $\alpha$ is not a constant
mapping. Let $\alpha=fa^m$ for some $f\in O_n$ and $0\le m\le n-1$. Then
$1f\ne nf$ for otherwise $\alpha$ would be constant, whence
$1\alpha\ne n\alpha$.

Conversely, let $\alpha\in\OP_n$. If $\alpha$ is a constant mapping then
$\alpha\in O_n$ and $\alpha=\alpha a^0\in O_n\cdot Z_n$. Now suppose that
$1\alpha\ne n\alpha$. We claim that $\alpha=fa^m$ for some $f\in O_n$ and
$0\le m\le n-1$. If $n\alpha>1\alpha$ then, by
Proposition~\ref{prop:2.3}, $\alpha$ is an order-preserving mapping and so
$\alpha=\alpha a^0$ is a required factorization. Suppose that
$n\alpha<1\alpha$ and write $1\alpha=k$. Put $m=n-k+1$. Then we define the
mapping $f$ as follows:
\[
xf=x\alpha+m\qquad(x\in X).
\]
Let $b\in X$. Hence $b(fa^{-m})=(bf)a^{-m}=bf-m=(b\alpha+m)-m=b\alpha$ and
so $fa^{-m}=\alpha$, which implies that $f=\alpha a^m\in\OP_n$. It remains
to prove that $f\in O_n$. Now $1f=1$ as $1f=k+m=1$, but
\[
nf=n\alpha+m\ne k+m,
\]
as $n\alpha\ne1\alpha=k$. Hence $1f<nf$, $f\in\OP_n$, and so $f\in O_n$,
by Proposition~\ref{prop:2.3}.
\end{proof}

From this result and the fact that $\OP_n=Z_n\cdot O_n$ we can conclude
that
\[
O_n\cdot(Z_n\cdot O_n)=O_n\cdot(\OP_n)=\OP_n
\]
and so $(O_n\cdot Z_n)^2=\OP_n$.

% ===============================================================
\section{Some properties of \texorpdfstring{$\OP_n$}{OPn}}
% ===============================================================

It is well known that $O_n$ is a regular semigroup. We can deduce from
this that $\OP_n$ is a regular semigroup.

Let $\alpha=a^kf\in\OP_n$ for some $f\in O_n$ and some $k$
$(0\le k\le n-1)$. It is routine to check that $f'a^{n-k}\in V(a^kf)$,
where $f'\in V(f)$ in $O_n$, though in Section~4 we show that not all
inverses of $a^kf$ are of this form. The next result follows from this
together with Theorem~\ref{thm:2.6}.

\begin{theorem}
$\OP_n$ is a regular semigroup. \qed
\end{theorem}

\begin{remark}
Since $\OP_n$ is regular, we have that
$\GR^{\OP_n}=\GR^{T_n}\cap(\OP_n\times\OP_n)$,
$\GL^{\OP_n}=\GL^{T_n}\cap(\OP_n\times\OP_n)$,
$\GH^{\OP_n}=\GH^{T_n}\cap(\OP_n\times\OP_n)$, that is, for
$\alpha,\beta\in\OP_n$, $\alpha\GR\beta\iff\ker\alpha=\ker\beta$,
$\alpha\GL\beta\iff X\alpha=X\beta$,
$\alpha\GH\beta\iff X\alpha=X\beta$ and $\ker\alpha=\ker\beta$, and since
$\OP_n$ is a finite semigroup we have that $\GD^{\OP_n}=\GJ^{\OP_n}$.
\end{remark}

Using the factorization of $\OP_n$ given by Theorem~\ref{thm:2.6} we can
see easily the following result relating the $\GL$-classes of $\OP_n$ and
$O_n$:

\begin{lemma}
Let $\alpha=a^kf$ and $\beta=a^mg$ for some $0\le k,m\le n-1$ and
$f,g\in O_n$. Then $\alpha\GL\beta$ in $\OP_n$ if and only if $f\GL g$ in
$O_n$.
\end{lemma}

However with $\GR$-classes the implication holds in neither direction for
$n>2$. For example, let $f=\tmap123223$ and $g=\tmap123113$. Then
$\ker f=\ker g$. However
$\alpha=a^2f=\tmap123322$, $\beta=ag=\tmap123131$ and
$\ker\alpha\ne\ker\beta$.

Next, with $f=\tmap123112$, $g=\tmap123133\in O_3$, let
$\alpha=af=\tmap123121$ and $\beta=a^2g=\tmap123313$. Then
$\ker\alpha=\ker\beta$, but $\ker f\ne\ker g$.

Let $\alpha,\beta\in\OP_n$. We have that $\alpha\GR\beta$ if and only if
$\ker\alpha=\ker\beta$, and kernel classes are convex, so there is a
one-to-one correspondence between $\GR$-classes of rank $t$, $(t\ge2)$,
and choices of a $t$-set from $X$ as follows: write the given $t$-set as
$\{a_1,a_2,\dots,a_t\}$ with $a_1<a_2<\cdots<a_t$; the kernel classes of
the corresponding $\GR$-class are then
$\{a_1,a_1+1,\dots,a_2-1\},\{a_2,a_2+1,\dots,a_3-1\},\dots,\{a_t,a_t+1,\dots,a_1-1\}$;
the $a_i$ are called the \emph{initial points} of the corresponding kernel
class; conversely, the common set of initial points of the kernel classes
of an $\GR$-class of $\OP_n$ yields a $t$-set of $X$, and these two
operations are mutually inverse.

We have $\alpha\GL\beta$ if and only if $X\alpha=X\beta$, so
$\GL$-classes also correspond to choices of a $t$-set from $X$, where
$t=|X\alpha|$. Hence for $t\ge2$ an $\GH$-class $H$ of $\OP_n$ can be coded
as a choice of an ordered pair of $t$-sets $(K,R)$ from $X$ where
$K=\{a_1,\dots,a_t\}$ is the common $t$-set of the initial points of the
kernel classes of the corresponding $\GR$-class and
$R=\{r_1,\dots,r_t\}$ $(r_1<r_2<\cdots<r_t)$ is the common range of all
mappings in the corresponding $\GL$-class.

Let $\alpha\in H$ with $a_i\alpha=r_1$, say. By Lemma~1.3~(b),
\[
(a_i,a_{i+1},\dots,a_t,a_1,\dots,a_{i-1})
\]
is cyclic, and since $\alpha\in\OP_n$ we have from Lemma~\ref{lem:2.1}
that the same is true of $(a_i\alpha,a_{i+1}\alpha,\dots,a_{i-1}\alpha)$.
However, $a_i\alpha=r_1=\min X\alpha$, and since the $a_j$ form a
transversal of the kernel classes of $\alpha$ we obtain:
\[
a_i\alpha<a_{i+1}\alpha<\cdots<a_t\alpha<a_1\alpha<\cdots<a_{i-1}\alpha.
\]
Hence $a_i\alpha=r_1$, $a_{i+1}\alpha=r_2,\dots,a_{i-1}\alpha=r_t$ and so
$\alpha\in H$ is defined by the pair of $t$-sets $(K,R)$ together with the
choice of $i$ such that $a_i\alpha=r_1$.

Conversely, for a given pair $(K,R)$ and a given choice of $i$ we can
construct $\alpha\in H$ as follows: write the kernel classes of $\alpha$
as $K_1,K_2,\dots,K_t$ where $a_j$ is the initial point of $K_j$. Define
$\alpha$ by $K_j\alpha=r_{j-i+1}$. In summary we have:

\begin{theorem}\label{thm:3.3}
For $t\ge2$ there is a one-to-one correspondence $\Phi$ between the set of
triples $(K,R,i)$ where $K$ and $R$ are ordered $t$-sets of $X$ and
$1\le i\le t$ and $\{\alpha\in\OP_n : |X\alpha|=t\}$ whereby
$(K,R,i)\longmapsto\alpha$, where each $a_j\in K$ is an initial point of a
kernel class of $\alpha$ and $a_j\alpha=r_{j-i+1}$. Moreover
$H_\alpha=\{\Phi(K,R,i) : i=1,\dots,t\}$, so $|H_\alpha|=t$. \qed
\end{theorem}

Suppose $\alpha,\beta\in\OP_n$ with $|X\alpha|=|X\beta|=t$ say. Let $K$ be
the set of initial points of the kernel classes of $\alpha$ and let
$R=X\beta$. Let $\gamma$ be any member of the $\GH$-class determined by the
pair $(K,R)$. Then $\alpha\GR\gamma\GL\beta$, so
$\alpha\GD^{\OP_n}\beta$. Conversely, if $\alpha\GD^{\OP_n}\beta$, then
$\alpha\GD^{T_n}\beta$ whence $|X\alpha|=|X\beta|$. Hence the next result
follows:

\begin{proposition}\label{prop:3.4}
Let $\alpha,\beta\in\OP_n$. Then $\alpha\GD^{\OP_n}\beta$ if and only if
$|X\alpha|=|X\beta|$. \qed
\end{proposition}

Next we give an alternative way of finding the order of $\OP_n$. Let
$\alpha\in\OP_n$ be such that $|X\alpha|=t$ for some $1\le t\le n$. Let
$D_t$ denote the $\GD$-class containing $\alpha$. By
Proposition~\ref{prop:3.4} we have $D_t=\{\beta\in\OP_n : |X\beta|=t\}$.
Hence, if $t=1$ we have $n$ constants mappings. By
Theorem~\ref{thm:3.3}, $\beta$ is determined by the pair
$(K_\beta,X\beta)$, where $K_\beta$ is the $t$-set of initial points of
the kernel classes of $\beta$, together with the choice of $i$ such that
$a_i\beta=r_1$ $(a_i\in K_\beta)$. Then for $t\ge2$ we have $C(n,t)$
possibilities for $K_\beta$, $C(n,t)$ choices for $X\beta$, and $t$
choices for $i$. Hence the order of $D_t$ $(t\ge2)$ is given by
\[
|D_t|=tC(n,t)C(n,t).
\]
Hence the order of $\OP_n$ is given by
\[
|\OP_n|=n+\sum_{t=2}^{n}tC^2(n,t).
\]
Note that
$\sum_{t=1}^{n}C(n-1,t-1)C(n,t)=\sum_{t=1}^{n}C(n-1,n-t)C(n,t)$ is the
coefficient of $x^n$ in the expansion of
$(1+x)^{n-1}(1+x)^n=(1+x)^{2n-1}$. Hence
$\sum_{t=1}^{n}C(n-1,t-1)C(n,t)=C(2n-1,n)$, and this leads to
$\sum_{t=1}^{n}tC^2(n,t)=nC(2n-1,n-1)$; thus we can recover
Corollary~\ref{cor:2.7}:
\[
|\OP_n|=n+nC(2n-1,n-1)-n^2.
\]

Next we give a description of the $\GH$-classes of $\alpha\in\OP_n$. Let
$\alpha\in E(\OP_n)$ and take $K_1,K_2,\dots,K_t$; $r_1,r_2,\dots,r_t$ and
$i$ as before. Since $\alpha\in E(\OP_n)$ the $\{r_j\}$ form a transversal
of the kernel classes, and so
\begin{equation}\label{eq:1}
r_{j-i+1}\in K_j,\qquad j=1,2,\dots,t.
\end{equation}
Let $\beta=\Phi(K,R,i-1)$. Hence,
\[
K_j\beta=a_j\beta=r_{j-(i-1)+1}=r_{j-i+2},\qquad j=1,2,\dots,t.
\]
Next we show by induction on $p$ that $a_{i-1}\beta^p=r_p$ for all
$p=1,2,\dots,t$. For $p=1$, we have that $a_{i-1}\beta=r_{i-1-i+2}=r_1$.
For $p\ge2$,
\[
a_{i-1}\beta^p=(a_{i-1}\beta^{p-1})\beta=r_{p-1}\beta.
\]
Put $j-i+1=p-1$ in \eqref{eq:1}; we get $j=p+i-2$. Hence
\[
a_{i-1}\beta^p=r_{p-1}\beta=K_{p+i-2}\beta=r_{p+i-2-i+2}=r_p.
\]
It follows that $|\gen{\beta}|\ge t$. Since $\beta\in H_\alpha$ and
$|H_\alpha|=t$ we see that $\gen{\beta}=H_\alpha$. Hence the next result
follows:

\begin{lemma}\label{lem:3.5}
Let $\alpha\in E(\OP_n)$. Let $H_\alpha$ be the $\GH$-class containing
$\alpha$. Then $H_\alpha=\gen{\beta}$ and $|H_\alpha|=|X\alpha|$. \qed
\end{lemma}

\begin{corollary}\label{cor:3.6}
The maximal subgroups of the $\GD$-class $D_t$ of $\OP_n$ are cyclic of
order $t$. In particular, every subgroup of $\OP_n$ is cyclic. \qed
\end{corollary}

% ===============================================================
\section{Conjugate subsemigroups of \texorpdfstring{$O_n$}{On}}
% ===============================================================

In this section we study some properties of subsemigroups of $O_n$
conjugate with respect to $Z_n$, and their links with the semigroup
$\OP_n$. Let
\[
O_n^k=\{a^{-k}fa^k : f\in O_n\},\qquad(0\le k\le n-1),
\]
where $a=(12\cdots n)$ is the fixed generator of $Z_n$.

Let $\alpha=a^mf\in\OP_n$ for some $f\in O_n$ and some $0\le m\le n-1$.
For each $k$, $(0\le k\le n-1)$, consider the mapping
$\varphi_k\colon O_n\longrightarrow\OP_n$ whereby
$h\varphi_k=a^{-k}ha^k$, $(h\in O_n)$. Let us write
$h^{a^k}=a^{-k}ha^k$. It is easy to verify that $\varphi_k$ is an
injective morphism with image $O_n^k$.

\begin{lemma}\label{lem:4.1}
$O_n^k$ is a submonoid of $\OP_n$ and $\varphi_k$ is an isomorphism
between $O_n$ and $O_n^k$.
\end{lemma}

\begin{proof}
This is an immediate consequence of the general fact that conjugation by
a unit is always a semigroup isomorphism.
\end{proof}

\begin{proposition}
$O_n^k$ is the monoid of all order-preserving mappings on $X$ with respect
to $\le_k$.
\end{proposition}

\begin{proof}
Let $h$ be any element of $O_n^k$. Hence $h=a^{-k}fa^k$ for some
$f\in O_n$. Let $i,j\in X$ such that $i\le_k j$. Hence $i-k\le j-k$, so
$ia^{-k}f=(i-k)f\le(j-k)f=ja^{-k}f$ since $f$ is order-preserving, and so
$ih=(i-k)f+k\le_k(j-k)f+k=jh$.

Conversely let $g$ be any order-preserving mapping on the chain
$(X,\le_k)$. Put $f=a^kga^{-k}$. By similar reasoning $f\in O_n$, so
$g=a^{-k}fa^k\in O_n^k$, and the proposition is proved.
\end{proof}

Let $T$ be any semigroup. Following standard terminology if $E=E(T)$ is
the set of idempotents of $T$, then $\gen{E}$ is called the \emph{core} of
$T$, and $T$ is said to be a \emph{semiband} if $T$ equals its own core.

We will denote by $G_b^S$ the $\GG$-class in $S$ containing $b$ where
$\GG$ is any Green's relations, $S$ is any semigroup and $b$ is an element
in $S$. We shall show that
$|E(R_f^{O_n})|=|E(R_{a^mf}^{\OP_n})|$. First, we need the following
result in order to show that $\varphi_k$ restricted to $E(R_f^{O_n})$ is a
bijection onto $E(R_{a^{-k}fa^k}^{\OP_n})$.

\begin{lemma}\label{lem:4.3}
Let $g\in O_n$, $e\in E(\OP_n)$ such that $e\GR^{\OP_n}g$. Then
$e\in O_n$ and $e\GR^{O_n}g$.
\end{lemma}

\begin{proof}
Observe that $\ker e=\ker g$, which is a convex partition of the chain
$X$. Let $i,j$ be any two elements on $X$ such that $i\le j$. If $i,j$ are
in the same kernel class then $ie=je$. Otherwise, $i\in K_i=xe^{-1}$ and
$j\in K_j=ye^{-1}$, for some distinct $x,y\in Xe$. Since $e$ is an
idempotent $x\in K_i$, $y\in K_j$ and since the kernel classes form a
convex partition of $X$, then $x\le y$, as required. The last assertion is
immediate since $O_n$ is a regular subsemigroup of $\OP_n$.
\end{proof}

\begin{lemma}\label{lem:4.4}
$\varphi_k$ restricted to $E(R_f^{O_n})$ is a bijection onto
$E(R_{a^{-k}fa^k}^{\OP_n})$.
\end{lemma}

\begin{proof}
Certainly $\varphi_k$ maps the first set injectively into the second. As
for surjectivity, take $e\in E(R_{a^{-k}fa^k}^{\OP_n})$. Then
$e^{a^{-k}}\varphi_k=e$, and $e^{a^{-k}}\in E(\OP_n)$. Now
$e^{a^{-k}}\GR^{\OP_n}a^ke\GR^{\OP_n}a^k(a^{-k}fa^k)=fa^k\GR^{\OP_n}f$.
But then $e^{a^{-k}}\in O_n$ and $e^{a^{-k}}\GR^{O_n}f$ by
Lemma~\ref{lem:4.3}, proving surjectivity.
\end{proof}

From this result we can now conclude with:

\begin{theorem}
Let $\alpha=a^kf\in\OP_n$ for some $f\in O_n$ and some $0\le k\le n-1$.
Then $|E(R_f^{O_n})|=|E(R_\alpha^{\OP_n})|$.
\end{theorem}

\begin{proof}
The result follows from Lemma~\ref{lem:4.4}, since
$\alpha\GR^{\OP_n}\alpha a^{-k}=a^kfa^{-k}$.
\end{proof}

Next we see that, in general,
$|E(L_f^{O_n})|\le|E(L_{a^mf}^{\OP_n})|$, since we have
$L_f^{O_n}\subseteq L_f^{\OP_n}=L_{a^mf}^{\OP_n}$. However we show that
the containment is strict for $n\ge3$. It is sufficient to observe this for
$n=3$. Take $\alpha\in\OP_3$ defined by $\alpha=a^0f$, $(f\in O_3)$ where
$f=\tmap123122$. Consider the mapping $\beta=\tmap123121\notin O_3$.
Clearly $\beta\in E(L_\alpha^{\OP_3})$ yet $\beta\notin E(L_f^{O_3})$.

Let $V^S(\alpha)$ be denote the set $V(\alpha)$ in the semigroup $S$. In
general we have that $|V^S(\alpha)|=|E(L_\alpha^S)|\,|E(R_\alpha^S)|$. In
our case, from the results above, we have
\[
|V^{O_n}(f)|=|E(L_f^{O_n})|\,|E(R_f^{O_n})|
<|E(L_{a^mf}^{\OP_n})|\,|E(R_{a^mf}^{\OP_n})|=|V^{\OP_n}(a^mf)|
\]
and so we can conclude that not all inverses of $a^mf$ are of the form
$f'a^{n-m}$ where $f'\in V^{O_n}(f)$.

It has long been known that the subsemigroup of all singular mappings of
$T_n$ is a semiband \cite{Howie66}, and it has been proved that $O_n$ is
also idempotent-generated \cite{Howie71}. Our next task is to find the core
of $\OP_n$. First we prove,

\begin{lemma}
Every $\alpha\in O_n$ has a fixed point.
\end{lemma}

\begin{proof}
Let $\alpha\in O_n$ and $A=\{x\in X : x\le x\alpha\}$. Note that $1\in A$
and so $A$ is nonempty. Let $a=\max A$. Hence $a\le a\alpha$ as $a\in A$,
and $a\alpha\le(a\alpha)\alpha$ as $\alpha\in O_n$. Thus $a\alpha\in A$
and so $a\alpha\le a$ by the maximality of $a$. Therefore $a=a\alpha$ as
required.
\end{proof}

Let $U=\bigcup_{k=0}^{n-1}O_n^k$ and consider $a^{-k}fa^k$ where
$f\in O_n$ and $0\le k\le n-1$. Now $x\in X$ is a fixed point of
$a^{-k}fa^k$ if and only if $(x-k)f=x-k$ which is equivalent to saying
that $x-k$ is a fixed point of $f$. Thus we have the next result:

\begin{lemma}\label{lem:4.7}
Let $\alpha=a^{-k}fa^k$ where $f\in O_n$ and $0\le k\le n-1$. Then
$x\in\F(\alpha)$ if and only if $x-k\in\F(f)$. In particular, if
$\alpha\in U$ then $\F(\alpha)\ne\emptyset$. \qed
\end{lemma}

As has often been observed each $\alpha\in T_n$ can be pictured as a
digraph on $n$ vertices with $ij$ an arc of $\alpha$ if $i\alpha=j$. Each
connected component of such a digraph is functional, meaning that it
consists of a unique cycle, together with a number of trees rooted around
the points of the cycle. Each component $C$ of the digraph of
$\alpha\in T_n$ has a \emph{cycle} $\mathrm{K}(C)$ defined by the property
that $x\in C$ is in $\mathrm{K}(C)$ if and only if there exists $p>0$ such
that $x\alpha^p=x$.

\begin{lemma}\label{lem:4.8}
Let $\alpha\in\OP_n$. Then the digraph of $\alpha$ cannot have a
non-trivial cycle and a fixed point.
\end{lemma}

\begin{proof}
Suppose that $\alpha$ has a fixed point $x$. Suppose $\alpha$ has a cycle
$C$ with $|C|=p\ge2$. Let $\le$ denote $\le_{x-1}$: note that $x$ is the
minimum of $(X,\le_{x-1})$. We claim that there exists $i\in C$ such that
$i<i\alpha$. Suppose not: then for each $i\in C$ we have:
\[
i>i\alpha>i\alpha^2>\cdots>i\alpha^{p-1}>i\alpha^p=i,
\]
which is a contradiction and so there is such an $i$. Since
$\alpha\in\OP_n$ and $x\in\F(\alpha)$, we obtain
$(x,i,i\alpha),(x,i\alpha,i\alpha^2),\dots,(x,i\alpha^j,i\alpha^{j+1})$
are cyclic. If for $j\ge0$ we had $i\alpha^j>i\alpha^{j+1}$ we have
$i\alpha^{j+1}<x$, which is a contradiction. Therefore
$i\alpha^j<i\alpha^{j+1}$, for all $j$; whence
$i<i\alpha<\cdots<i\alpha^{p-1}<i\alpha^p=i$, which is a contradiction.
Hence the digraph of $\alpha$ cannot have a proper cycle and a fixed
point.
\end{proof}

\begin{theorem}
Let $\alpha\in\OP_n$. Then the digraph of $\alpha$ cannot have two cycles
of different length.
\end{theorem}

\begin{proof}
Suppose $\alpha$ has a $t$-cycle $C_t$ and an $r$-cycle $C_r$ with $t<r$.
Consider $\alpha^t\in\OP_n$. Each point of $C_t$ is fixed by $\alpha^t$
and so $\F(\alpha^t)\ne\emptyset$. However, for $x\in C_r$,
$x\alpha^t\ne x$, but $x$ is a point in a cycle of $\alpha^t$ since
$x\alpha^{tr}=x(\alpha^r)^t=x$, which contradicts Lemma~\ref{lem:4.8}.
\end{proof}

\begin{lemma}\label{lem:4.10}
Let $\alpha\in\OP_n$. If $\F(\alpha)\ne\emptyset$ then the digraph of
$\alpha$ is a forest and each component $C$ associated with a fixed point
of $\alpha$ is an interval.
\end{lemma}

\begin{proof}
If $\F(\alpha)\ne\emptyset$, then by Lemma~\ref{lem:4.8} the digraph of
$\alpha$ cannot have a non-trivial cycle. Thus each component has a
trivial cycle and so the digraph of $\alpha$ is a forest. Now let $C$ be a
component with $C=C(x)$, $x$ the fixed point in $C$. Take $k$ sufficiently
large so that $C(x)=\{y\in X : y\alpha^k=x\}$. Then
$C(x)=x(\alpha^k)^{-1}$; $\alpha^k\in\OP_n$ and so $x(\alpha^k)^{-1}$ is
an interval by the Convexity Lemma.
\end{proof}

\begin{lemma}\label{lem:4.11}
Let $\alpha\in\OP_n$ such that $\F(\alpha)\ne\emptyset$. Let $C$ be any
component of $\alpha$. Then $C\alpha\subseteq C$, there exist $i,j\in X$
such that $C=[i,j]$, and $\alpha$ restricted to $C$ is an order-preserving
with respect to $\le_{i-1}$.
\end{lemma}

\begin{proof}
It follows immediately from the definition of component that
$C\alpha\subseteq C$, and from Lemma~\ref{lem:4.10} that $\alpha$
restricted to $C$ has a fixed point $x$. Suppose first that $C=X$. Then
$C=[x,x-1]$. Suppose that $a\le_{x-1}b$. If $b\alpha=x$ then
$x=x\alpha=a\alpha=b\alpha$, by the Convexity Lemma. If
$x<_{x-1}b\alpha$ then $(x,a\alpha,b\alpha)$ is cyclic, by
Lemma~\ref{lem:2.1}, so $x\le_{x-1}a\alpha\le_{x-1}b\alpha$. This proves
that $\alpha$ is order-preserving with respect to $\le_{x-1}$.

Suppose now that $C\ne X$ and choose any $z\in X\setminus C$. Observe, by
definition of component, that $z\alpha\notin C$. By
Lemma~\ref{lem:4.10}, $C=[i,j]$ for some $i,j\in X$. In particular note
that $j<_{i-1}z$ and $j<_{i-1}z\alpha$. Let $y\in C$. If $y\le_{i-1}x$
then $(y,x,z)$ is cyclic, so $(y\alpha,x,z\alpha)$ is cyclic, by
Lemma~\ref{lem:2.1}, so $y\alpha\le_{i-1}x$. Similarly if $x\le_{i-1}y$
then $x\le_{i-1}y\alpha$. Now suppose that $a,b\in C$ and $a<_{i-1}b$. If
$a\le_{i-1}x\le_{i-1}b$ then, from above,
$a\alpha\le_{i-1}x\le_{i-1}b\alpha$. If $a,b\le_{i-1}x$ then $(a,b,x)$ is
cyclic so that $(a\alpha,b\alpha,x)$ is cyclic and, from above,
$a\alpha,b\alpha\le_{i-1}x$. If $a\alpha=x$ then $a\alpha=x\alpha$ and so
$a\alpha=b\alpha=x$ by the Convexity Lemma. Otherwise $a\alpha<_{i-1}x$,
$b\alpha\le_{i-1}x$, and so by cyclicity, $a\alpha\le_{i-1}b\alpha$.
Similarly, if $x\le_{i-1}a,b$ then $a\alpha\le_{i-1}b\alpha$, completing
the proof that $\alpha|_C$ is order-preserving with respect to
$\le_{i-1}$.
\end{proof}

\begin{corollary}\label{cor:4.12}
Let $\alpha\in\OP_n$ such that $\F(\alpha)\ne\emptyset$. Then
$\alpha\in U$.
\end{corollary}

\begin{proof}
Let $C$ be any component of $\alpha$, so $C=[i,j]$ for some $i,j$ and
$\alpha|_C\colon C\longrightarrow C$ is order-preserving with respect to
$\le_{i-1}$, by Lemma~\ref{lem:4.11}. Let $D$ be any other component (if
it exists). Then, by Lemma~\ref{lem:4.11} again, $D=[k,l]$ for some $k,l$
and $\alpha|_D\colon D\longrightarrow D$ is order-preserving with respect
to $\le_{k-1}$. But
\[
\le_{i-1}\cap\,(D\times D)\;=\;\le_{k-1}\cap\,(D\times D),
\]
since $i\notin D$ and $D$ is an interval, so that $\alpha|_D$ is
order-preserving with respect to $\le_{i-1}$. Thus $\alpha$ is
order-preserving with respect to $\le_{i-1}$, so lies in $U$.
\end{proof}

The previous lemmas and the description of the fixed points of
$a^{-k}fa^k$ in terms of those of $f$ now yield the next result:

\begin{theorem}\label{thm:4.13}
$U=\{\alpha\in\OP_n : \F(\alpha)\ne\emptyset\}$.
\end{theorem}

\begin{proof}
If $\alpha\in U$ then $\F(\alpha)\ne\emptyset$ by Lemma~\ref{lem:4.7}.
The converse follows from Corollary~\ref{cor:4.12}.
\end{proof}

$U$ does not contain $a=(12\cdots n)$ nor its non-idempotent powers.
Moreover there are elements in $\OP_n\setminus Z_n$ which are not in $U$
either. For example, let
$\alpha=\begin{pmatrix}1&2&3&4\\2&4&4&2\end{pmatrix}\in\OP_4$. We have
$\alpha\notin Z_4$ as $\alpha$ is not a power of $a=(1234)$, the fixed
generator of $Z_4$; however $\F(\alpha)=\emptyset$.

\begin{corollary}
$E(\OP_n)\subseteq U$.
\end{corollary}

\begin{proof}
If $\alpha\in E(\OP_n)$ then $\F(\alpha)\ne\emptyset$, and so by
Theorem~\ref{thm:4.13}, $\alpha\in U$.
\end{proof}

From this result, if $\alpha\in\gen{E}$ then $\alpha\in\gen{U}$. However
$O_n$ is a semiband \cite{Howie71} which implies that $O_n^k$ is a
semiband for all $k=0,1,\dots,n-1$ by Lemma~\ref{lem:4.1}. Let
$\alpha\in U$, so that there exists $k$ $(0\le k\le n-1)$ such that
$\alpha\in O_n^k$. Since $O_n^k$ is a semiband, then $\alpha$ is a product
of elements in $E(O_n^k)$ and hence $\alpha$ is a product of elements in
$E(\OP_n)$. Therefore $\alpha\in\gen{E}$. Thus
$\gen{U}\subseteq\gen{E}$ as $U\subseteq\gen{E}$. Hence the next result
follows:

\begin{corollary}
Let $E=E(\OP_n)$ and $U=\bigcup_{k=0}^{n-1}O_n^k$. Then
$\gen{E}=\gen{U}$. \qed
\end{corollary}

Next we describe $\gen{U}$, the idempotent-generated subsemigroup of
$\OP_n$.

\begin{lemma}\label{lem:4.16}
Let $\alpha\in\OP_n$. Then $\alpha a^p\in U$ for some $p$
$(0\le p\le n-1)$.
\end{lemma}

\begin{proof}
Let $\alpha\in\OP_n$. Using Theorem~\ref{thm:2.6}, $\alpha=a^mg$ for some
$0\le m\le n-1$ and $g\in O_n$. Then
$\alpha a^{n-m}=a^mga^{n-m}=a^{-(n-m)}ga^{n-m}\in O_n^{n-m}\subseteq U$.
\end{proof}

\begin{lemma}\label{lem:4.17}
Let $R\subseteq X$. Then the mapping $f$ whereby:
\[
xf=\begin{cases} xa & \text{if } x\in R\\ x & \text{if } x\notin R\end{cases}
\]
is in $\OP_n$. Moreover if $R$ is a proper subset, $f\in U$.
\end{lemma}

\begin{proof}
Note that for any $i$ such that $1\le i\le n-1$, we have
$(i-1)f\le i\le if$; therefore $1f\le 2f\le\cdots\le(n-1)f$. Now $nf=n$ or
$nf=1$; in the first case $1f\le 2f\le\cdots\le(n-1)f\le nf$ and
$f\in O_n$; in the second case $nf\le 1f\le 2f\le\cdots\le(n-1)f$ and
$f\in\OP_n$. In either case $f\in\OP_n$. Now if $R$ is a proper subset
then $\F(f)\ne\emptyset$ as $\F(f)=X\setminus R$ and so by
Theorem~\ref{thm:4.13}, $f\in U$.
\end{proof}

\begin{lemma}\label{lem:4.18}
Suppose that $u\in\gen{U}$, $u\ne1$. Then $ua\in\gen{U}$.
\end{lemma}

\begin{proof}
Note that the identity mapping is the only member of $U$ of full rank as
$U=\bigcup_{m=0}^{n-1}a^{-m}O_na^m$ and the identity mapping is the only
member of $O_n$ of full rank. Hence $Xu\ne X$. Take $R=Xu$ in
Lemma~\ref{lem:4.17}. We claim that $ua=uf$, from which the result follows
from Lemma~\ref{lem:4.17}. Clearly $f|_{Xu}=a|_{Xu}$. Let $x$ be any
element of $X$. Then we have that $x(ua)=(xu)a=(xu)f=x(uf)$ and so
$ua=uf$.
\end{proof}

Finally we have a description of the core of $\OP_n$.

\begin{theorem}
$\gen{U}=\gen{E(\OP_n)}=(\OP_n\setminus Z_n)\cup\{1_X\}$.
\end{theorem}

\begin{proof}
Let $\alpha\in(\OP_n\setminus Z_n)\cup\{1_X\}$. Clearly
$1_X\in\gen{U}$. Suppose $\alpha\ne1_X$. By Lemma~\ref{lem:4.16}, there
exists $p$ such that $\alpha a^p\in\gen{U}$, and clearly
$\alpha a^p\ne1_X$. Take $m$ such that $p+m\equiv0\pmod n$. Using
Lemma~\ref{lem:4.18} we therefore obtain
$\alpha=\alpha a^{p+m}=\alpha a^pa^m\in\gen{U}$.

Conversely, $U$, and so $\gen{U}$, contains no members of $\OP_n$ of rank
$n$ (except $1_X$), and so contains no other mappings of full rank; in
particular
\[
Z_n\cap\gen{U}=\{1_X\},
\]
as required.
\end{proof}

Consequently the core of $\OP_n$ is $\gen{U}$ and we can conclude that
$\OP_n$ is not a semiband, although, like $T_n$, the submonoid of singular
mappings of $\OP_n$ is idempotent-generated.

There are mappings in the core of $\OP_n$ without fixed point; for example
the mapping exhibited before Corollary~4.14. Note that the digraph of that
mapping has one component. For an example whose digraph has two
components, take
$\alpha=\begin{pmatrix}1&2&3&4&5&6\\3&4&5&2&3&3\end{pmatrix}\in\OP_6$. We
have $\alpha\notin U$ as $\F(\alpha)=\emptyset$ but $\alpha\in\gen{U}$.

% ===============================================================
\section{Orientation-reversing mappings}
% ===============================================================

Let $\alpha\in T_n$. We say that $\alpha$ is an
\emph{orientation-reversing} mapping on $X$ if the sequence
$(1\alpha,2\alpha,\dots,n\alpha)$ is anti-cyclic. The collection of all
orientation-reversing mappings on $X$ will be denoted by $\OR_n$.

Let $\gamma$ be the reflection whereby $i\longmapsto n+1-i$ $(i\in X)$.
Note that $\gamma$ is an involution on $X$. Also, $\gamma\in\OR_n$ as
$(1\gamma,2\gamma,\dots,n\gamma)=(n,n-1,\dots,1)$ is anti-cyclic.

Let $\Gamma_r$ [resp.\ $\Gamma_l$]$\colon T_n\longrightarrow T_n$ be
defined by $\alpha\longmapsto\alpha\gamma$ [resp.\
$\alpha\longmapsto\gamma\alpha$] $(\alpha\in T_n)$. The following result is
easily proved.

\begin{lemma}
$\Gamma_r,\Gamma_l$ are involutions on $T_n$ which map $\OP_n$ onto
$\OR_n$ and $\OR_n$ onto $\OP_n$.
\end{lemma}

\begin{corollary}
$\OP_n\cdot\OR_n=\OR_n\cdot\OP_n=\OR_n$ and
\[
(\OR_n)^2=(\OP_n)^2=\OP_n.
\]
\end{corollary}

Let $P_n$ be denote the collection of all mappings in $\OP_n\cup\OR_n$. By
Corollary~5.2, we have that:

\begin{theorem}
$P_n=\OP_n\cup\OR_n$ is a submonoid of $T_n$. \qed
\end{theorem}

From Lemma~1.1 we obtain the following

\begin{lemma}\label{lem:5.4}
$\OP_n\cap\OR_n=\{\alpha\in\OP_n : |X\alpha|\le2\}$. \qed
\end{lemma}

From this we can conclude that

\begin{theorem}
$|P_n|=2|\OP_n|-n-2C^2(n,2)$.
\end{theorem}

\begin{proof}
Note that $|P_n|=|\OP_n|+|\OR_n|-|\OP_n\cap\OR_n|$. From
Lemma~\ref{lem:5.4} and the discussion following
Proposition~\ref{prop:3.4}, $|\OP_n\cap\OR_n|=n+2C^2(n,2)$ and since
$|\OP_n|=|\OR_n|$ the result follows.
\end{proof}

\begin{remark}
A simple calculation now yields
\[
|P_n|=\frac{(2n)!}{n!\,(n-1)!}-\frac{n}{2}\bigl(n^3-2n^2+5n-2\bigr)
\]
and this gives $P_n=T_n$ for $n\le3$. This may be a reason why many
interesting counterexamples of semigroups of transformations require 4
base points.
\end{remark}

Let $\beta\in\OR_n$ and write $\beta$ as $\alpha\gamma$ where
$\alpha\in\OP_n$. It is routine to check that
$\gamma\alpha'\in V(\alpha\gamma)$ in $P_n$, where $\alpha'\in V(\alpha)$
in $\OP_n$. (Conversely, given $\beta'\in V(\beta)$ in $P_n$ then it is
straightforward to check that $\beta'=\alpha'\gamma$ for some
$\alpha'\in V(\alpha)$ in $\OP_n$.) Hence the next result follows:

\begin{theorem}
$P_n$ is a regular monoid. \qed
\end{theorem}

Consequently, from this result the usual description of Green's relations
$\GL$, $\GR$ and $\GH$ applies again. (See Remark that follows
Theorem~3.1.)

\begin{lemma}
Let $\alpha,\beta\in P_n$. Then $\alpha\GD^{P_n}\beta$ if and only if
$|X\alpha|=|X\beta|$.
\end{lemma}

\begin{proof}
If $\alpha\GD^{P_n}\beta$, then $\alpha\GD^{T_n}\beta$ whence
$|X\alpha|=|X\beta|$.

Conversely, suppose that $|X\alpha|=|X\beta|$. If $\alpha,\beta\in\OP_n$
then $\alpha\GD^{\OP_n}\beta$ by Proposition~\ref{prop:3.4}. If
$\alpha,\beta\in\OR_n$ then $|X\alpha\gamma|=|X\beta\gamma|$, so that, by
Proposition~\ref{prop:3.4},
\[
\alpha\GR^{P_n}\alpha\gamma\GD^{\OP_n}\beta\gamma\GR^{P_n}\beta,
\]
so certainly $\alpha\GD^{P_n}\beta$. Similarly, if only one of
$\alpha,\beta$ lies in $\OP_n$, say without loss of generality
$\alpha\in\OP_n$, $\beta\in\OR_n$, then
$\alpha\GD^{\OP_n}\beta\gamma\GR^{P_n}\beta$. In each case
$\alpha\GD^{P_n}\beta$.
\end{proof}

\begin{lemma}
Let $\alpha\in P_n$ be such that $|X\alpha|=t$ for some $3\le t\le n$. Let
$H_\alpha$ be the $\GH$-class containing $\alpha$. Then $|H_\alpha|=2t$.
\end{lemma}

\begin{proof}
Let $\alpha\in P_n$ be such that $|X\alpha|=t$ for some $1\le t\le n$. Let
$D_t^{P_n}$ denote the $\GD$-class containing $\alpha$. By the previous
lemma, $D_t^{P_n}=\{\beta\in P_n : |X\beta|=t\}$. By
Lemma~\ref{lem:5.4}, for $t\le2$ we have
$D_t^{P_n}\cap\OP_n=D_t^{P_n}\cap\OR_n$. For $t\ge3$,
$D_t^{P_n}=(D_t^{P_n}\cap\OP_n)\cup(D_t^{P_n}\cap\OR_n)$. (The union is
disjoint.) Since $\Gamma_r$ restricted to $D_t^{P_n}\cap\OP_n$ is a
bijection onto $D_t^{P_n}\cap\OR_n$ then
$|D_t^{P_n}\cap\OP_n|=|D_t^{P_n}\cap\OR_n|$. Since
$D_t^{P_n}\cap\OP_n=D_t^{\OP_n}$ then $|D_t^{P_n}|=2|D_t^{\OP_n}|$. In
$P_n$, $D_t^{P_n}$ has the same number of $\GR$-, $\GL$-classes as
$D_t^{\OP_n}$ and so the number of $\GH$-classes in $D_t^{P_n}$ is equal
the number of $\GH$-classes in $D_t^{\OP_n}$. Since all $\GH$-classes
within the one $\GD$-class have the same cardinality, then the order of
every $\GH$-class in $D_t^{P_n}$ is twice the order of every $\GH$-class
in $D_t^{\OP_n}$, whence the result follows from
Corollary~\ref{cor:3.6}.
\end{proof}

\begin{theorem}
For $t\ge3$ the maximal subgroups of $D_t^{P_n}$ are dihedral groups of
order $2t$.
\end{theorem}

\begin{proof}
Recall that the dihedral group $D_n$ of degree $n$ is a subgroup of $S_n$
$(n\ge3)$ of order $2n$ and generators $x$ and $y$ $(x,y\in D_n)$
satisfying the following conditions: $x^n=1_X$, $y^2=1_X$, $x^k\ne1_X$ if
$0<k<n$ and $yx=x^{-1}y$; conversely any group satisfying these conditions
is isomorphic to $D_n$. (See, for example, \cite{Hungerford}.)

Let $\alpha\in E(P_n)$ be such that $|X\alpha|=t$ for some $t$. Take
$K_1,K_2,\dots,K_t$ as the kernel classes of $\alpha$ and
$r_1,r_2,\dots,r_t$ as the range of $\alpha$ as before such that
$K_j\alpha=r_j$, $j=1,2,\dots,t$. Next we describe the maximal subgroups of
$D_t^{P_n}$ and we show that they are dihedral groups of order $2t$ with
generators $\delta$ and $\beta$ $(\delta,\beta\in P_n)$ satisfying
$\beta^t=\alpha$, $\delta^2=\alpha$, $\beta^k\ne\alpha$ if $0<k<t$ and
$\delta\beta=\beta^{-1}\delta$, where $\beta^{-1}$ is the group inverse of
$\beta$ in $H_\alpha$. For this consider the mapping $\delta$ defined as
follows:
\[
K_j\delta=r_{t-j+1},\qquad j=1,2,\dots,t.
\]
Clearly $X\delta=X\alpha$ and $\ker\delta=\ker\alpha$. Moreover
$\delta\in\OR_n$ as the sequence $(r_t,r_{t-1},\dots,r_1)$ is
anti-cyclic. Also, consider the mapping $\beta$ whereby
\[
K_j\beta=r_{j+1},\qquad j=1,2,\dots,t.
\]
Note that $\beta$ is the mapping used in the proof of
Lemma~\ref{lem:3.5} and $\beta\in\OP_n$, $\beta^t=\alpha$,
$\beta^k\ne\alpha$ if $0<k<t$. Also for all $j=1,2,\dots,t$ we have
\[
K_j\delta^2=r_{t-j+1}\delta=r_{t-(t-j+1)+1}=r_j,
\]
and so $\delta^2=\alpha$. Now the group inverse of $\beta$ in $H_\alpha$
is $\beta^{-1}$ defined by $K_j\beta^{-1}=r_{j-1}$, $j=1,2,\dots,t$.
Finally, $(K_j\beta^{-1})\delta=r_{j-1}\delta=r_{t-(j-1)+1}$. On the other
hand we obtain
\[
(K_j\delta)\beta=r_{t-j+1}\beta=r_{(t-j+1)+1}
\]
and so $\delta\beta=\beta^{-1}\delta$.
\end{proof}

\subsection*{Acknowledgements}
The first author gratefully acknowledges support from the Foundation
Calouste Gulbenkian and the University of Tr\'as-os-Montes e Alto Douro.
Also, we would like to acknowledge a number of helpful detailed
improvements suggested by the referee.

% ===============================================================

\end{document}